\documentclass{elsart3-1}

\usepackage{amsmath}
\usepackage{amssymb}
\usepackage{color}

\usepackage[english,francais]{babel}

\newtheorem{theorem}{Theorem}[section]

\newtheorem{e-proposition}[theorem]{Proposition}

\newtheorem{e-definition}[theorem]{Definition\rm}
\newtheorem{remark}{\it Remark\/}

\renewcommand{\theequation}{\arabic{equation}}
\newcommand{\dt}[1]{\frac{\partial #1}{\partial t}}
\newcommand{\dx}[1]{\frac{\partial #1}{\partial x}}
\newcommand{\dxx}[1]{\frac{\partial^2 #1}{\partial x^2}}

\newcommand{\erf}{\mathrm{erf}}

\def\og{\leavevmode\raise.3ex\hbox{$\scriptscriptstyle\langle\!\langle$~}}
\def\fg{\leavevmode\raise.3ex\hbox{~$\!\scriptscriptstyle\,\rangle\!\rangle$}}

\journal{the Acad\'emie des sciences}
\begin{document}
\centerline{}
\begin{frontmatter}


\selectlanguage{english}
\title{On a price formation free boundary model by Lasry \& Lions}


\selectlanguage{english}
\author[authorlabel1]{Luis A. Caffarelli}
\ead{caffarel@math.utexas.edu}
\author[authorlabel2,authorlabel3]{Peter A. Markowich}
\ead{P.A.Markowich@damtp.cam.ac.uk}
\author[authorlabel2]{Jan-F. Pietschmann}
\ead{J.Pietschmann@damtp.cam.ac.uk}

\address[authorlabel1]{Department of Mathematics, Institute for Computational Engineering and Sciences, University of Texas at Austin, USA}
\address[authorlabel2]{DAMTP, University of Cambridge, Cambridge CB3 0WA, UK}
\address[authorlabel3]{Faculty of Mathematics, University of Vienna, 1090 Vienna, Austria}


\medskip
\begin{center}
{\small Received *****; accepted after revision +++++\\
Presented by }
\end{center}

\begin{abstract}
\selectlanguage{english}
We discuss global existence and asymptotic behaviour of a price formation free boundary model introduced by Lasry \& Lions in 2007. Our results are based on a construction which transforms the problem into the heat equation with specially prepared initial datum. The key point is that the free boundary present in the original problem becomes the zero level set of this solution. Using the properties of the heat operator we can show global existence, regularity and asymptotic results of the free boundary.

\vskip 0.5\baselineskip

\selectlanguage{francais}
\noindent{\bf R\'esum\'e} \vskip 0.5\baselineskip \noindent
{\bf A propos d'un mod\`ele de formation des prix de Lasry \& Lions}
Nous discutons l'existence globale et le comportement asymptotique d'un mod\`ele de formation des prix \`a fronti\`ere libre introduit par Lasry \& Lions en 2007. Nos r\'esultats sont bas\'es sur une transformation qui relie ce probl\`eme \`a l'equation de la chaleur pour une condition initiale bien choisie. L'id\'ee principale est que la fronti\`ere libre pr\'esente dans le mod\`ele original est donn\'ee par la ligne de niveau z\'ero de cette solution. En utilisant les propri\'et\'es de l'op\'erateur de la chaleur, nous pouvons montrer des r\'esultats de r\'egularit\'e sur la fronti\`ere libre, et donc l'existence d'une solution globale.

\end{abstract}
\end{frontmatter}


\selectlanguage{english}
\section{Introduction}
\label{sec:intro}
This paper is concerned with a mean field game model in economics which was introduced in a paper by J.-M. Lasry and P.-L. Lions, cf. \cite{Lasry2007}. The setup consists of a (large) group of buyers and a (large) group of vendors trading a good at a certain price $p(t)$, with a fixed transaction cost $a$. The model is given by a non-linear parabolic free boundary evolution equation that describes the dynamical behaviour of the densities of buyers and vendors which in turn define the price. It is set up on the whole real line, i.e. the price can, in principle, take arbitrarily large values. The model is given by the equation
 \begin{align}\label{eq:f_inv}
 \dt{f} - \dxx{f} &= \lambda(t)(\delta(x-p(t)+a) - \delta(x-p(t)-a)),\, x\in\mathbb{R},t\in\mathbb{R}_+\\
\lambda(t) &= -\dx{f}(p(t),t),\; f(p(t),t)=0,\\\label{eq:initial}
 f(x,0) &= f_I,\, p(0) = p_0,\;\text{for some $p_0$ in $\mathbb{R}$},
\end{align}
with compatibility conditions at time $t=0$:
\begin{align*}
 \text{(A1)}\qquad\qquad f_I(p_0) = 0\text{ and }f_I(x) > 0\text{ for }x < p_0 \text{ and }f_I(x) < 0\text{ for }x > p_0.
\end{align*}
For the following we assume that $f_I$ is in $L^1(\mathbb{R})$ and bounded. This model has been studied in a number of papers, cf. \cite{Chayes2009}, \cite{Ganzalez2011} and \cite{Markowich2009}. Here we shall present the first global existence result of a smooth solution on the whole real line. In the following, we shall denote by $f = f^+ - f^-$ the decomposition of a function into its positive and negative part.

\section{Connection to the Heat Equation}
In this section we shall prove that there is a one to one relation between solutions of the aforementioned FBP \eqref{eq:f_inv}-\eqref{eq:initial} and solutions of the heat equation supplemented with specially prepared initial data. This will lead to a global existence result in a very elegant way.
\begin{theorem} Let $f=f(x,t)$ be a solution solution of \eqref{eq:f_inv}-\eqref{eq:initial} on the time interval $[0,T]$ with $T>0$. Then there exists a linear transformation from $f$ to a function $F=F(x,t)$, being a solution of the heat equation, such that the graph of the zero level set of $F$ is $p(t)$. By reversing the transformation, each solution of the heat equation such that the zero level set of the solution is a smooth graph for $0\le t \le T$ can be transformed into a solution of the FBP with the same $p(t)$.
\end{theorem}
The construction is based on the observation that the second derivative of $-f^-$ at the free boundary $p(t)$ is precisely the negative value of the weighted delta mass centered at $p(t)+a$, as it appears in the equation \eqref{eq:f_inv}. Analoguosly, the second derivative of $f^+$ is the negative of the weighted delta mass of the equation \eqref{eq:f_inv}, centered at $p(t)-a$.\\
{\it Proof.} Let $f_I=f_I(x)$ be a given initial datum satisfying assumption (A1). Let $f=f(x,t)$ be the solution of \eqref{eq:f_inv}-\eqref{eq:initial} in the time interval $[0,T]$ (such a solution exists due to \cite[Theorem 2.6]{Markowich2009}). Now we define
\begin{align*}
 F(x,t) = \left\{\begin{array}{cc}
                  \phantom{-}\sum_{n=0}^\infty f^+(x+na,t), & x < p(t),\\
		  -\sum_{n=0}^\infty f^-(x-na,t), & x > p(t).
                  \end{array}\right.
\end{align*}
We remark that due to the boundedness of $f$ these sums converge in $\mathcal{D}'(\mathbb{R}\times [0,\infty))$. It is very easy to check that $F$ satisfies, in the sense of distributions, the heat equation with initial datum $F(x,t=0)=:F_I(x)$, given by \eqref{eq:Fi}. Clearly, the free boundary $p=p(t)$ is now the zero level set of $F$. Now consider a given $F=F(x,t)$, solution of the heat equation in $[0,T]$. Assume the initial datum is of the form
\begin{align}\label{eq:Fi}
 F_I(x) = \left\{\begin{array}{cc}
                  \phantom{-}\sum_{n=0}^\infty f_I^+(x+na), & x < p_0,\\
		  -\sum_{n=0}^\infty f_I^-(x-na), & x > p_0,
                  \end{array}\right.
\end{align}
for an arbitrary function $f_I$ satisfying (A1). Then, we can construct a solution of the FBP \eqref{eq:f_inv}-\eqref{eq:initial} with the initial datum $f_I$ in the following way:
\begin{align*}
 f(x,t) = \left\{\begin{array}{cc}
                  F^+(x,t) - F^+(x+a), & x < p(t),\\
		  -F^-(x,t) + F^-(x-a), & x > p(t).
                  \end{array}\right.
\end{align*}
Again, by construction, the zero level set of $F$ becomes the free boundary of \eqref{eq:f_inv}-\eqref{eq:initial}.
\begin{theorem}[Global Existence] There exists a unique smooth solution $f=f(x,t)$ of \eqref{eq:f_inv}-\eqref{eq:initial} for $t\in [0,\infty)$. Furthermore, $p\in\mathcal{C}([0,\infty))$. 
\end{theorem}
{\it Proof.} Let $F_I$ be the transformed initial datum corresponding of $f_I$ and let $F$ be the solution of the heat equation with initial datum $F_I$. Abusing notation, we denote by $p=p(t)$ the zero level set of $F$. First we note that oscillations of $p(t)$ yielding a 'fat' free boundary cannot occur as they contradict the $x$-analyticity of solutions of the heat equation. Furthermore, due to \cite[Lemma 2.9]{Markowich2009} we know that $f_x(p(t),t) < 0$ for all $t>0$ (by the Hopf Lemma) and the min-max principle implies that $p=p(t)$ is the graph of a function. Hence we only need to exclude the existence of $t^*$ such that $|p(t)|$ becomes unbounded as $t\rightarrow t^*$. We write
\begin{align}\label{eq:duhamel}
 F(x,t) &= \int_{-\infty}^\infty G(t,z)F_I(x-z)\;dz
=\int_{-\infty}^{x-p_0} G(t,z)F_I^-(x-z)\;dz - \int_{x-p_0}^\infty G(t,z)F_I^+(x-z)\;dz,
\end{align}
where 
 $G(t,x) = \frac{1}{\sqrt{4\pi t}}\exp{\left(\frac{x^2}{4t}\right)}$
is the 1-d heat kernel. Due to boundedness of $f$ and its construction, $F$ grows at most linearly at $|x|=\infty$. Thus second term on the right hand side in \eqref{eq:duhamel} tends to zero as $x\rightarrow +\infty$. For the first term we have
\begin{align*}
\int_{-\infty}^{x-p_0} G(t,z)F_I^-(x-z)\;dz=\int_{p_0-x}^{\infty} G(t,z)F_I(x+z)\;dz\ge C\int_{p_0}^{p_0+a} |F_I(x+z)|\;dz.
\end{align*}
Due to \eqref{eq:Fi} we have
\begin{align*}
\int_{p_0}^{p_0+a} |F_I(x+z)|\;dz= \sum_{n=0}^\infty \int_{x+p_0-na}^{x+p_0-(n-1)a} f_I^-(y)\;dy \ge \mathrm{const} > 0.
\end{align*}
Thus for $x$ large enough, this term dominates in \eqref{eq:duhamel} and thus $F(\cdot,t)$ becomes negative. By the same argument we show that for large negative $x$, $F(\cdot,t)$ becomes positive and thus there must exist a unique $x$ with $-\infty < x < \infty$ such that $F(x,t)=0$. From these arguments we conclude that $p(t)$ is defined and continuous for all $t$.
\begin{remark} 
A similar analysis produces solutions of the Neumann problem in the interval $[-L,L]$ and certain examples of non-existence. In this case the associated solution of the heat equation satisfies the unusual Neuman type boundary condition  $F_x(\pm L,t)=F_x(\pm L \mp a,t)$.
\end{remark}

\section{Asymptotic Behaviour}
From now on we assume an initial datum $f_I$ with $p_0=0$ and $a=1$. For the following we define the function $\erf(u) := \frac{1}{\sqrt{4\pi}}\int_u^\infty e^{-\frac{x^2}{4}}\;dx$. In this section, we shall prove:
\begin{theorem} Let $f=f(x,t)$ be a solution of \eqref{eq:f_inv}-\eqref{eq:initial}. If $M^+ := \int_{-\infty}^0 f^+(z)\;dx \neq \int_0^\infty f^-(z)\;dz =: M^-$, then $p(t) \sim \sqrt{t}q_\infty$ with $\erf(q_\infty) = \left.M^-\right/ M^+$ as $t\rightarrow\infty$. If $M^- = M^+$, i.e. the total mass of $f$ is zero, then
\begin{align*}
 p(t) = \frac{\int_{-\infty}^\infty z|f(z)|\;dz}{M^+ + M^-} + O\left(\frac{1}{\sqrt{t}}\right).
\end{align*}
\end{theorem}
{\it Proof.} From \eqref{eq:Fi} and \eqref{eq:duhamel} we obtain at $x=p(t)$
\begin{align*}
 0&= -\int_0^\infty\frac{1}{\sqrt{4\pi t}}\sum_{n=0}^\infty\exp{\left(-\frac{|p(t)-z+n|^2}{4t}\right)}f^-(z)\;dz + \int_{-\infty}^0\frac{1}{\sqrt{4\pi t}}\sum_{n=0}^\infty\exp{\left(-\frac{|p(t)-z+n|^2}{4t}\right)}f^+(z)\;dz. 
\end{align*}
The sums in the above equation for the free boundary can be interpreted as Riemann sums converging to integrals. We easily obtain:
\begin{align*}
 \int_0^\infty\erf \left(\frac{p(t)}{\sqrt{t}} - \frac{z}{\sqrt{t}} - \frac{1}{2\sqrt{t}}\right)f^-(z)\;dz = \int_{-\infty}^0\left(1-\erf\left(\frac{p(t)}{\sqrt{t}}- \frac{z}{\sqrt{t}}+\frac{1}{2\sqrt{t}}\right)\right)f^+(z)\;dz + O\left(\frac{1}{t}\right).
\end{align*}
Now we define $q(t):= p(t)/\sqrt{t}$ and conclude the proof of the first part of the theorem. We remark that in particular $q_\infty \neq 0$ if $M^+ \neq M^-$. To prove the second part, we note that as $t\rightarrow\infty$ we have
\begin{align*}
 \mathrm{erf}\left(\frac{p(t)-z-\frac{1}{2}}{\sqrt{t}}\right)\sim \frac{1}{2} - \frac{1}{\sqrt{4\pi t}}\left(p(t)-z-\frac{1}{2}\right) + O\left(\frac{1}{t}\right).
\end{align*}
Using the assumption $M^+ = M^-$, this allows us to write \eqref{eq:duhamel} as
\begin{align*}
 0=F(p(t),t) = \frac{p(t)}{\sqrt{4\pi t}}\left(\int_{-\infty}^\infty \left|f(z)\right|\;dx - \left(\int_{0}^\infty zf^-(z)\;dz + \int_{-\infty}^0 zf^+(z)\;dz\right)\right) + O\left(\frac{1}{t}\right).
\end{align*}
for $t$ large enough. Thus we immediately obtain
\begin{align*}
 p(t) = \frac{\int_{-\infty}^\infty z\left|f(z)\right|\;dz}{M^++M^-} + O\left(\frac{1}{\sqrt{t}}\right),
\end{align*}
which concludes the proof.
\section*{Acknowledgements}
We acknowledge supported by Award No. KUK-I1-007-43, made by King Abdullah University of Science and Technology (KAUST), by the Leverhulme Trust via {\em Kinetic and mean field partial differential models for socio-economic processes\/} (PI Peter Markowich) and by the Royal Society through the Wolfson Research Merit Award of Peter Markowich. Luis Caffarelli acknowledges supprt from the Division of Mathematical Sciences of the NSF.




\end{document}